\input amstex
\documentstyle{amsppt}
\def\br{\Bbb R}
\def\bc{\Bbb C}
\def\bn{\Bbb N}
\def\bq{\Bbb Q}
\def\bz{\Bbb Z}
\def\cd{\Cal D}
\def\cf{\Cal F}
\def\ch{\Cal H}
\def\cm{\Cal M}
\def\cR{\Cal R}
\def\cl{\Cal L}
\def\clo{\Cal L^2(\Omega)}
\def\clm{\Cal L^2(\mu)}

\topmatter
\title Estimates on the Spectrum of Fractals Arising From Affine Iterations \endtitle
\rightheadtext{A Spectrum for Affine Fractals}
\author Palle E.T. Jorgensen and Steen Pedersen \endauthor
\address Department of Mathematics, University of Iowa, Iowa City, IA 52242, USA \endaddress
\email jorgen\@math.uiowa.edu \endemail
\address Department of Mathematics, Wright State University, Dayton, OH  45435 \endaddress
\email spedersen\@desire.wright.edu \endemail
\thanks Work supported in part by the U.S. National Science Foundation \endthanks
\abstract In the first section we review recent results on the harmonic analysis of fractals generated by {\it iterated function systems\/} with emphasis on {\it spectral duality}. Classical harmonic analysis is typically based on groups whereas the fractals are most often not groups. We show that nonetheless those fractals that come from iteration of affine mappings in $\br^d$ have a spectral duality which is instead based on approximation and a certain dual affine system on the Fourier transform side. The present work is based on iteration of {\it frame estimates\/} (which have been studied earlier for regions in $\br^d$). Our emphasis is on new results regarding the interplay between the limit-fractal on the one hand, and on the other the corresponding regions in $\br^d$ which generate iterated function systems of contractive affine mappings. As an application of our frame results, we obtain a classification of a certain type of spectral pairs.
\endabstract
\endtopmatter

\document
\head 1. Introduction \endhead
Let $\Omega$ be a measurable subset of $\br^d$ $(d=1,2,\ldots)$ of finite positive measure, i.e., for fixed $d$, $m_d$ denotes the $d$-dimensional Lebesgue measure and we assume, $0<m_d(\Omega) <\infty$. Then the Hilbert space  $\clo$ is formed relative to the usual inner product
$$
\langle f,g\rangle _\Omega := \int_\Omega \overline{f(x)}\,g(x)\, dx \tag 1.1
$$
where $dx:= dx_1 \ldots dx_d$, and $\|f \|_\Omega^2 := \langle f,f \rangle_\Omega$. If $\lambda =(\lambda_1, \ldots , \lambda_d) \in \br^d$, then we shall need the exponential functions
$$
e_\lambda(x) := e^{i2\pi\lambda \cdot x}
=e^{i2\pi (\lambda_1 x_1 + \cdots +\lambda_d x_d)}. \tag 1.2
$$
These functions are defined on all of $\br^d$ of course, but, when $\Omega$ is given as above, they {\it restrict\/} and define elements in $\clo$, viz., $\chi_\Omega e_\lambda$ where $\chi_\Omega$ denotes the indicator function of $\Omega$. It will be convenient to denote these restricted functions also by $e_\lambda$. If $\Lambda$ is a subset of $\br^d$ such that the corresponding exponential functions $\{e_\lambda : \lambda\in \Lambda\}$ form an orthogonal basis for the Hilbert space $\clo$, then we say that $(\Omega,\Lambda)$ is a {\it spectral pair}. Such pairs have been studied extensively recently, see e.g., \cite{Fu}, \cite {Jo1}, \cite{JoPe3--4}, and \cite{LaWa}. Defining the Fourier transform
$$
\hat f(\lambda):= \int_\Omega \overline{e_\lambda(x)}\,f(x)\, dx
= \langle e_\lambda, f \rangle_\Omega \tag 1.3
$$
we then get
$$
\sum_{\lambda\in\Lambda} |\hat f(\lambda)|^2 = m_d(\Omega) \int _\Omega |f(x)|^2 \, dx
=m_d(\Omega)\|f\|_\Omega^2.
\tag 1.4
$$
Of course, the set $\Omega$ is implicit in the notation on the left-hand side of \rom{(1.4)}, as $\Omega$ enters into the definition of the transform, $f \mapsto \hat f$, see \rom{(1.3)}.

It turns out that the structure of {\it spectral pairs\/} is surprisingly rigid, and it is also related to {\it tiling\/} questions, see \cite{Fu}, \cite{JoPe3} and \cite{LaWa}. The {\it generalized spectral pairs\/} (see details below) turn out to form a much richer class. They arise from a relaxation of \rom{(1.4)}, transforming it instead into a {\it frame estimate\/}, see also \cite{Ben1--2, BeTo, Ga, DuSc, Beu}. We say that $(\Omega, \Lambda)$ is a {\it generalized spectral pair\/} (g.s.p.) in $\br^d$ if $\Omega$ is as specified, and if there are constants, $k$, $K$, positive and finite, but depending on both $\Omega$ and $\Lambda$, such that
$$
k\|f\|_\Omega^2 \leq \sum_{\lambda\in\Lambda} |\hat f (\lambda)|^2
\leq K\|f\|_\Omega^2. \tag 1.5
$$
It can be checked that if $\Lambda$ is a set satisfying \rom{(1.5)} for $\Omega$ as specified, then $\Lambda$ must be discrete. 

\proclaim{Lemma 1.1} Let $(\Omega,\Lambda)$ be a generalized spectral pair, and let $$
K=\Lambda^\circ :=\{ t\in \br^d : t\cdot \lambda \in \bz,\ \forall \lambda\in\Lambda\}.
$$
Then $K$ is discrete. \endproclaim

\demo{Proof} The set $\{e_\lambda : \lambda\in \Lambda\}$ is total in $\clo$. Hence by \cite {Fu}, $\Lambda$ cannot be contained in any hyperplane in $\br^d$; in particular, there exists a basis $\{ v_1, \ldots, v_d\}$ for $\br^d$ with each $v_j \in \Lambda$. Let $w_1,\ldots,w_d$ be a dual basis, viz., $w_j \cdot v_i =\delta_{ji}$, for $i$, $j=1,\ldots,d$. Then we get the containment,
$$\align
K&\subseteq \{ t\in \br^d : t\cdot v_j \in \bz,\  j=1,\ldots,d\}  \\
 & = \left\{ \sum_{j=1}^d t_j w_j : t_j\in \bz\right\},
\endalign$$
with the latter set a lattice in $\br^d$; hence $K$ is discrete. \qed\enddemo

Note, we only need $\{e_\lambda : \lambda \in \Lambda\}$ to be total in $\clo$ for $\Lambda^\circ$ to be discrete---not any of the other ingredients in the definition of a generalized spectral pair are in fact needed for this part.

It is also know that, for the spectral pairs, if $0\in \Lambda$, then the group
$$
\Lambda^\circ := \{s\in \br^d : s\cdot \lambda \in \bz,\ \forall \lambda \in \Lambda \}
\tag 1.6
$$
is of rank $d$; in other words, $\Lambda^\circ$ is a {\it lattice\/} in $\br^d$. If $\Lambda$ {\it is also a lattice}, we have just the familiar situation of multivariable Fourier transforms; but there are many extremely interesting cases where $(\Omega,\Lambda)$ is a spectral pair, but $\Lambda$ is {\it not\/} a lattice. We take up this question in Sections 2.3--2.4 below which provides in fact a classification of some of the possibilities. Although the results there are about spectral pairs, the arguments use g.s.p.'s and iteration.

If $(\Omega, \Lambda)$ is a generalized spectral pair, then there is a well defined operator $F:= F_{(\Omega,\Lambda)}$ from $\clo$ into $\ell^2(\Lambda)$, given by
$$
(Ff)(\lambda):= \hat f(\lambda),\qquad \lambda\in\Lambda \tag 1.7
$$
and it is called the associated {\it frame operator}. We say that the system is {\it exact\/} if $F_{(\Omega,\Lambda)}$ maps {\it onto\/} $\ell^2(\Lambda)$, and that it is {\it tight\/} if $k=K$ in \rom{(1.5)}. It follows from Hilbert space theory that the system is {\it exact\/} if and only if there are no non-zero sequence solutions $(\xi_\lambda)_{\lambda\in\Lambda}$ in $\ell^2(\Lambda)$ such that the corresponding function $f(x)=\sum_{\lambda \in \Lambda} \xi_\lambda e_\lambda (x)$ is zero as an element in $\clo$.

Let $R$ be a real $d$ by $d$ matrix which is {\it expansive}; i.e., we assume that the eigenvalues $\lambda_i$ of $R$ all satisfy $|\lambda_i|>1$. If $R=(r_{ij})$, the corresponding transpose $(r_{ji})$ will be denoted $R^*$. Two systems of affine maps will be considered as follows: Let $B$ and $L$ be given finite subsets in $\br^d$, and set
$$\align
\sigma_b(x) &:= R^{-1}(x+b) =R^{-1}x + R^{-1}b,\tag 1.8\\
\tau_\ell(s) &:= R^*(s+\ell) =R^*s + R^*\ell,\tag 1.9 \\
\sigma (\Omega) &:= \bigcup_{b\in B} \sigma_b(\Omega),\tag 1.10\\
\intertext{and}
\tau(\Lambda) &:=\bigcup _{\ell \in L} \tau_\ell (\Lambda).\tag 1.11
\endalign$$
The affine mappings \rom{(1.8)--(1.9)} are defined for $x$, $s\in \br^d$, and they are indexed by the respective (finite) sets $B$ and $L$. In \rom{(1.8)--(1.11)}, we are applying these maps to subsets of $\br^d$, e.g., the components of some pair $(\Omega, \Lambda)$ as above, or to more general sets. We shall be primarily interested in the case when the unions in \rom{(1.10)--(1.11)} are non-overlapping, i.e., when
$$
\sigma_b(\Omega)\cap \sigma_{b'}(\Omega)= \emptyset\qquad\text{for all $b\neq b'$ in $B$}\tag 1.12
$$
and
$$
\tau_\ell (\Lambda)\cap \tau_{\ell '}(\Lambda)= \emptyset
\qquad\text{for all $\ell \neq \ell '$ in $L$}. \tag 1.13
$$
As far as \rom{(1.12)} is concerned, we will relax it and allow overlap on sets of $d$-dimensional Lebesgue measure zero, i.e.
$$
m_d(\sigma_b(\Omega)\cap \sigma_{b'}(\Omega))=0
\qquad\text{for all $b\neq b'$ in $B$}.\tag 1.14
$$
But it turns out that the harmonic analysis of systems subjected to \rom{(1.14)} may be reduced to that of an associated system which in fact satisfies the stronger property \rom{(1.12)}, so we will restrict here to the latter case. When \rom{(1.14)} holds we say that the  triple $(R,\Omega,B)$ is an affine system with {\it no overlap}, and if \rom{(1.13)} holds that $(R^*, \Lambda,L)$ is a {\it dual affine system\/} with non-overlap. Starting with \rom{(1.10)--(1.11)}, we define recursively, $\Omega_0 := \Omega$, $\Lambda_0 :=\Lambda$,
$$
\Omega_{n+1} := \sigma(\Omega_n),\tag 1.15
$$
and
$$
\Lambda_{n+1} := \tau(\Lambda_n).\tag 1.16
$$
For each $n$, let $\mu_n$ be the normalized measure coming from restricting Lebesgue measure $m_d$ to $\Omega_n$, i.e.,
$$
\mu_n := \chi_{\Omega_n} \times m_d (\Omega_n)^{-1} \times m_d. \tag 1.17
$$
The interpretation of \rom{(1.17)} is
$$
\mu_n(\Delta)=m_d(\Omega_n)^{-1} m_d(\Delta \cap \Omega_n)\qquad\text{for all Borel sets $\Delta$,}
$$
or equivalently,
$$
\int f\,d\mu_n = m_d(\Omega_n)^{-1} \int f\chi_{\Omega_n}\,dm_d\qquad\text{for all Borel functions $f$.}
$$
It can be checked that then 
$$
\mu_{n+1} = \frac{1}{|B|} \sum_{b\in B} \mu_n \circ \sigma_b^{-1} \tag 1.18
$$
where $|B|$ denotes the number of elements in the (finite) set $B$, and $\mu_n \circ \sigma_b^{-1}$ is the measure given by,
$$
\mu_n \circ \sigma_b^{-1}(E) := \mu_n (\sigma _b^{-1}(E)),
$$
where
$$
\sigma_b^{-1}(E) = \{ x: \sigma_b(x)\in E\},
$$
or, equivalently,
$$
\int f\left( d\mu_n \circ \sigma_b^{-1}\right) = \int (f\circ \sigma_b)\, d\mu_n \tag 1.19
$$
for Borel subsets $E$, and Borel functions $f$. Using then a standard theorem (see \cite{Hut} and \cite{Fal}) there is a unique limit measure $\mu$ (i.e., $\mu_n \rightarrow \mu$ in the Hutchinson-Hausdorff metric), and
$$
\mu = \frac{1}{|B|} \sum_{b\in B} \mu \circ \sigma_b^{-1}.\tag 1.20
$$
This is typically a fractal measure, and the object of the present paper is the harmonic analysis of the corresponding Hilbert space $\clm$, i.e., the usual completion, now defined from
$$
\| f\|_\mu ^2 := \int |f|^2\,d\mu .\tag 1.21
$$

The simplest examples of this would be the following Cantor constructions: 
$$
d=1,\qquad R=4,\qquad \Omega=I =[0,1]\qquad\text{the unit-interval, and}\qquad B=\{0,2\};
$$
or, alternatively, 
$$
R=3\qquad\text{and}\qquad B=\{0,2\}.
$$
The corresponding fractal dimensions are $1/2=\ln 2/ \ln 4$, respectively, $0.631\ldots =\ln2 / \ln3$. As we show in \cite{JoPe5--7}, there are sets $\Lambda\subset \br$ such that for the first example, $\clm$ has a harmonic analysis which arises as a limit of the {\it spectral pairs\/} indicated in the construction \rom{(1.15)--(1.16)} above, but this construction fails for the second example (the ternary Cantor set), and we will see in this paper that there is a limit construction which is still analogous, but it is based instead on the {\it generalized spectral pairs}, see \rom{(1.5)} above.

Following \rom{(1.3)}, we may define a transform $\hat \mu$ for Borel probability measures $\mu$ on $\br^d$:
$$
\hat\mu (\lambda) := \int \overline{e_\lambda(x)}\, d\mu(x)
$$
and
$$
\widehat{f\, d\mu}(\lambda) :=\int \bar e_\lambda f\, d\mu \tag 1.22
$$
as functions in $\br^d$. In his recent papers \cite{Str1--2}, Strichartz considers integrals on the form
$$
\frac{1}{m_d(E_j)^\alpha} \int_{E_j} \left|\widehat{f\, d\mu}(s)\right|^2 \, ds \tag 1.23
$$
where $\alpha$ is a number (typically a {\it fraction\/} related to the {\it fractional dimension\/} of $\mu$), $E_j$ a family of sets, $E_j \subset E_{j+1}$, such that $0<m_d(E_j) < \infty$, and $\bigcup_j E_j =\br^d$. He shows that $\alpha$ may be chosen such that if $\underline{\cl}$ and $\overline{\cl}$ denote the respective $\liminf_j$ and $\limsup_j$ in \rom{(1.23)}, then there are positive constants $k$, $K$ such that
$$
k\| f\|_\mu^2 \,\leq\, \underline{\cl}\, \leq \,\overline{\cl}\, \leq\, K \|f\|_\mu^2. \tag 1.24
$$
Compare this to the {\it frame\/} property \rom{(1.5)} above. We shall be interested here in {\it discrete versions\/} of these spectral estimates; specifically, we shall give conditions on the dual affine system, as in \rom{(1.8)--(1.9)}, which guarantees an estimate as that of Strichartz in \rom{(1.23)}, but now with the $\liminf$ and $\limsup$'s defined instead from the sequence
$$
\left( \frac{|B|}{|L|}\right)^j \sum_{\lambda\in\Lambda_j} \left| \widehat{(f\,d\mu_j)}(\lambda)\right|^2.\tag 1.25
$$
The transform under the summation is defined in \rom{(1.22)} and is the same as the one used by Strichartz. But neither there, nor here, is it clear {\it when\/} the $j\rightarrow \infty$ limit exists. In our setting, we illustrate this issue with results and examples in Section 3 below.

Frame estimates as in \rom{(1.5)} play a role in the theory of wavelets, and fractal geometry, which is analogous to that played by the classical G\aa{}rding inequalities in the theory of elliptic partial differential equations (p.d.e.). A main result of ours (details in Sections 2 and 3 below) states that, in the present setting, there is a finite-dimensional matrix $A=\cm^*\cm$ such that a given positive lower spectral bound for $A$ implies an {\it a priori\/} ``smoothness'' for our fractal measure $\mu$, and the corresponding sequence of transforms \rom{(1.22)} which is the counterpart of a classical p.d.e.-G\aa{}rding estimate. Specifically we show that a lower bound on the spectrum of $A$ implies the following {\it a priori\/} estimate (for some $\epsilon\in \br_+$)
$$
\liminf_{j\rightarrow\infty} \left(\frac{|B|}{|L|}\right)^j
\sum_{\lambda\in\Lambda_j} \left|\widehat{(f\, d\mu_j)}\,(\lambda)\right|^2
\geq \epsilon \int |f|^2\, d\mu \tag 1.26
$$
for our fractal measure $\mu$ ($= \lim_{j \rightarrow\infty}\mu_j$ in the Hutchinson metric), which in turn arises as a solution to
$$
\mu=\frac{1}{|B|} \sum_{b\in B} \mu \circ \sigma_b^{-1} \tag 1.27
$$
within the probability measures on $\br^d$. The equation \rom{(1.27)}  plays a role here which is quite analogous to that of an elliptic p.d.e.\ in the classical case. This is perhaps more clear if \rom{(1.27)} is viewed as a special case of a more general class of {\it difference equations}, viz.,
$$
\mu=\sum_{b\in B} w_b (x)\,\mu\circ\sigma_b^{-1} \tag 1.28
$$
where $\{ w_b (\cdot)\}_{b\in B}$ is a given family of functions on $\br^d$. For more on this, see also \cite{Str2}. The interpretation of \rom{(1.28)} is the identity
$$
\int f\,d\mu= \sum_{b\in B} \int w_b (\sigma_b(x))f(\sigma_b(x))\,d\mu(x) 
\qquad\text{(for $\forall f\in C_c(\br^d)$)}.
$$

Our {\it dual affine systems\/} are said to be {\it elliptic\/} if the two finite translation sets $B$ and $L$ determine a matrix $A=\cm^*\cm$ with (positive) spectrum bounded away from zero. Our analogy to the G\aa{}rding case is the result (with suitable technical conditions) that ellipticity implies the {\it a priori\/} estimate \rom{(1.26)}.

\head 2. The Initial Step for the Iteration \endhead
\subhead 2.1 Pairs in $\br^d$ \endsubhead
Let $d$ be fixed, $d=1,2,\ldots$, and let $m_d$ be Lebesgue measure. In Section 1 we introduced pairs of subsets $(\Omega, \Lambda)$ of $\br^d$ such that $\Omega \subset\br^d$ is measurable and $0<m_d(\Omega)<\infty$, and the set $\Lambda$ serves as a {\it spectrum\/} for $\Omega$ in one of the two senses that we described; see especially \rom{(1.5)} which corresponds to $(\Omega,\Lambda)$ being a {\it generalized spectral pair}. If $d=1$, then, in applications, $\Omega$ may be {\it time\/} and $\Lambda$ {\it frequency}. From quantum mechanics, we may have $d>1$, and points $x=(x_1,\ldots,x_d)$ in $\Omega$ describe the {\it position\/} of a particle or system with many degrees of freedom, whereas points $\lambda=(\lambda_1, \ldots,\lambda_d)$ in $\Lambda$ describe the corresponding dual vector of {\it momentum\/} variables.

In \rom{(1.10)} above, we considered the step from one set $\Omega$ to a second set, $\Omega + B$, which is the union of {\it translates\/} of $\Omega$ itself. The simplest special case of this would be when $B=\{ 0,a\}$, and the new set in $x$-space is, $\Omega \cup (\Omega+a)$, assuming again that
$$
\Omega \cap (\Omega+ a) =\emptyset.\tag 2.1
$$
Here $\Omega + a=\{x+a : x\in \Omega\}$.

We shall be interested in generating {\it new generalized spectral pairs\/} from old ones, and eventually in the iteration of the process, leading thereby to affine {\it fractals}. The aim is to study the interplay between the original configuration in $\br^d$ on the one hand, and the iteration limit on the other hand. The original configuration has a harmonic analysis based on $\br^d$ as an additive group and Lebesgue measure as Haar measure. The iteration fractal will have the canonical measure $\mu$ which is our substitute for Haar measure.

One reason for the {\it generalized\/} spectral pairs is that, even if the starting point is a {\it spectral pair}, then the first, or one of the later affine iteration steps, typically takes us outside the restricted class of spectral pairs, but still we stay within the {\it generalized\/} class (i.e., g.s.p.). The generalized spectral pairs $(\Omega, \Lambda)$ are determined by the corresponding frame operator, $F=F_{(\Omega,\Lambda)}$ from \rom{(1.3)} and \rom{(1.7)}, i.e., the operator, $F: f\mapsto (\hat f(\lambda))_{\lambda\in \Lambda}$ of $\clo$ into $\ell^2(\Lambda)$; and we will now study how this operator $F$ changes when we pass from $\Omega$ to 
$$
\Omega_1:= \Omega \cup (\Omega +a)\tag 2.2
$$
where $a\in \br^d$ is fixed, but chosen subject to \rom{(2.1)}. Recall the operator $F$ maps from $\clo$ into $\ell^2(\Lambda)$. It is bounded, and it has zero kernel, i.e., the only solution to, $Ff=0$, $f\in \clo$, is $f=0$. It follows then from Hilbert space theory \cite{BeFr} that $F^*F$ is a {\it self-adjoint\/} operator in $\clo$. It has a spectral resolution
$$
F^*F = \int_{\br^+} \xi E(d\xi)\tag 2.3
$$
where $E(\cdot)$ is an orthogonal projection-valued measure on $\br^+$. We will study the above mentioned problem by checking how \rom{(2.3)} changes when we pass from $\Omega$ to $\Omega_1$, where again $\Omega_1$ is determined as in \rom{(2.2)} above. But we have in mind also more general steps, $\Omega\mapsto \Omega_1$, and we will follow through and supplying an accompanying modification, $\Lambda \mapsto \Lambda_1$, i.e., we will identify $\Lambda_1$ such that $(\Omega_1,\Lambda_1)$ is again a {\it generalized spectral pair\/} with an associated selfadjoint frame operator $F_1^*F_1$, now in $\cl^2(\Omega_1)$ relative to Lebesgue measure on $\Omega_1$. The simplest special case of this construction is (see \cite{JoPe7}) the case $d=1$, $\Omega=[0,1]=I$, $a\in \br$, $a>1$, and $\Lambda= \bz+\beta$, for some fixed $\beta \in \br$. It is clear that then $(\Omega,\Lambda)$ is a spectral pair. This is just the classical Fourier transform.

\example{Example 2.1.1} We considered in \cite{JoPe7} $a>1$ for the non-overlap property \rom{(2.1)}, and we defined $\Omega_1$ and $\Lambda_1$ as follows:
$$
\Omega_1 := I\cup (I+a) 
 =[0,1]\cup [a,a+1]\tag 2.4
$$
and
$$
\Lambda_1 :=(\bz + \beta)\cup a^{-1}\bz
$$
and we showed that $(\Omega_1,\Lambda_1)$ is a generalized spectral pair iff $a\in \bq$ and $\beta \notin \bz +a^{-1}\bz$. The operator $F_1^*F_1$ is still well defined also if $a$ is {\it irrational}, but we showed that the lower bound in the spectrum is then $0$, and therefore that the lower bound in the frame estimate \rom{(1.5)} is then violated for $\Omega_1 = [0,1] \cup [a,a+1]$, {\it when $a$ is irrational}. These facts will follow from our Theorem 2.1.2 below, which is the corresponding general $\br^d$ case. We will also give a completely explicit formula for the spectral resolution of the modified operator $F_1^*F_1$.

If $(\Omega,\Lambda)$ is a given generalized spectral pair (g.s.p.\ for short) in $\br^d$, then it is clear from \rom{(1.5)} that, for every vector $\ell \in \br^d$, the pair $(\Omega,\Lambda+\ell)$ is also a g.s.p.; and, similarly, we may make a translation of $\Omega$ and not affect the g.s.p.\ property. It follows that the condition, $0\in \Lambda$, may be added, and the g.s.p.\ property is unaffected. In the case $0\in \Lambda$, we shall need the dual discrete set $\Lambda^\circ$ from \cite{JoPe3} and Lemma 1.1 above. We shall also need the similar duality \rom{(2.5)} below, but now defined over the {\it rationals}, rather than the integers, i.e.,
$$
\Lambda^\circ(\bq) :=\{ s\in \br^d : s\cdot \lambda \in \bq,\ \forall \lambda\in \Lambda\}. \tag 2.5
$$
\endexample

\proclaim{Theorem 2.1.2} Let $(\Omega,\Lambda)$ and $(\Omega_2,\Lambda_2)$ be two given spectral pairs in $\br^d$ and assume that $0$ is in both $\Lambda$ and $\Lambda_2$ and $\Omega \subset \Omega_2$. Consider points $a$ and $\beta$ in $\br^d$ such that $a$ satisfies the disjointness property \rom{(2.1)}, and now form a third pair $(\Omega_1,\Lambda_1)$ as follows:
$$
\Omega_1 := \Omega\cup (\Omega +a) \qquad\text{and}\qquad
\Lambda_1 := (\Lambda +\beta)\cup \Lambda _2. \tag 2.6
$$
If $a\in \Lambda_2^\circ$, then the new pair $(\Omega_1,\Lambda_1)$ is a g.s.p.\ if and only if
\roster
\item"{(i)}" $a\notin (\Lambda +\beta)^\circ$  and
\item"{(ii)}" $a\in \Lambda^\circ (\bq)$.\endroster
In any case the spectrum of $F_1^*F_1$ is the set $\{r_\pm(\lambda) : \lambda\in \Lambda\}$ where
$$
r_\pm (\lambda) = 2+\alpha \pm \sqrt{(2+\alpha)^2 -4\alpha (1-\cos 2\pi(\beta+\lambda)\cdot a)} \tag 2.7
$$
and
$$
\alpha = 2 m_d (\Omega_2)/m_d(\Omega). \tag 2.8
$$
\endproclaim

\demo{Proof} Let $\tau$ be the period-2 automorphism of $\Omega_1$ which sends $x\in \Omega$ to $x+a$, and $x\in \Omega+a$ to $x-a$; and define
$$
Tf := f\circ \tau\qquad\text{for all } f\in \cl^2(\Omega_1).\tag 2.9
$$
This operator is unitary in $\cl^2(\Omega_1)$ with spectral subspaces
$$\ch_\pm := \{f\in \cl^2 (\Omega_1) : Tf = \pm f\} \tag 2.10
$$
and 
$$
\ch_+ \oplus \ch_- = \cl^2 (\Omega_1). \tag 2.11
$$
For elements in $f$ in $\cl^2(\Omega_1)$ we may consider the transform $\hat f$ defined in \rom{(1.3)}, and the integration being over $\Omega_1$, but we shall also need the restriction of $f$ to $\Omega \subset \Omega_1$ and the corresponding transform, now with integration only over $\Omega$. The latter will be denoted $\tilde f(\cdot)$, specifically
$$
\tilde f (\lambda) := \int_\Omega \overline{e_\lambda(x)}\,f(x)\, dx,\qquad\text{for $\lambda \in \br^d$}.\tag 2.12
$$
We shall be using the following shorthand notation
$$
\langle \lambda,x\rangle := e_\lambda (x) 
= e^{i2\pi \lambda \cdot x} \tag 2.13
$$
defined for all $\lambda$, $x \in \br^d$. When $f_\pm \in \ch_\pm$ (respectively) then the $\hat f_\pm$ transforms compute with integration over $\Omega_1 = \Omega \cup (\Omega +a)$ as follows
$$
\hat  f_+(\lambda) = (1+ \overline{\langle \lambda,a\rangle} ) \tilde f_+(\lambda) \tag 2.14
$$
and
$$
\hat f_- (\lambda) = (1-\overline{\langle \lambda, a\rangle}) \tilde f_- (\lambda) \tag 2.15
$$
for all $\lambda \in \br^d$. Since $a\in \Lambda_2^\circ$ we conclude that
$$
\hat f_- (\lambda_2) =0\qquad\text{for all $\lambda_2\in \Lambda_2$}.
$$
From (i) in the theorem, we conclude that the two terms in $\Lambda_1 = (\Lambda + \beta) \cup \Lambda_2$ are disjoint, so the $\Lambda_1$-summation in
$$
\sum_{\lambda_1 \in \Lambda_1} \left| \hat f(\lambda_1)\right|^2 \tag 2.16
$$
may be calculated separately for $\lambda \in \Lambda +\beta$, and then for $\lambda\in \Lambda_2$. If the two sets in the union $\Lambda_1:=(\Lambda+\beta) \cup \Lambda_2$ are not disjoint, we will still get a frame estimate, but then in the following modified form,
$$
\frac12 s \leq\sum_{\lambda\in \Lambda_1} |\hat f(\lambda)|^2 \leq s
$$
where $s:= \sum_{\lambda\in \Lambda+\beta} |\hat f(\lambda)|^2 + \sum_{\lambda\in\Lambda_2} |\hat f(\lambda)|^2$; but the disjointness is implied by \rom{(i)} and $a\in \Lambda_2^\circ$. With the $\ch_\pm$-decomposition, we get
$$
\hat f (\lambda) = (1+\overline{\langle \lambda,a\rangle})\tilde f_+(\lambda)
  +(1-\overline{\langle \lambda, a \rangle})\tilde f_-(\lambda)
\qquad\text{for all }\lambda\in \Lambda + \beta; \tag 2.17
$$
and
$$
\hat f(\lambda_2) = 2\tilde f_+(\lambda_2)\qquad\text{for all }\lambda_2\in\Lambda_2. \tag 2.18
$$

It follows that
$$\align
\sum_{\lambda_2\in \Lambda_2} \left| \hat f(\lambda_2)\right|^2
  & = 4 \sum_{\lambda_1 \in \Lambda_2} \left| \tilde f_+ (\lambda_2)\right|^2\\
  & =4m_d (\Omega_2) \times \int_{\Omega_2} \left| \left. f_+ \right|_\Omega \right|^2\, dx\\
  &= 4m_d(\Omega_2) \times \int_\Omega |f_+ |^2 \, dx\\
  & =4\frac{m_d(\Omega_2)}{m_d(\Omega)} \sum_{\lambda \in \Lambda +\beta}
         \left| \tilde f_+(\lambda) \right|^2
\endalign$$
where we used that $\Omega\subset\Omega_2$ and the assumption that both $(\Omega, \Lambda)$ and $(\Omega_2,\Lambda_2)$ are given to be spectral pairs. But if they are known instead to be only g.s.p.'s, then the same calculation would give the estimate
$$
\sum_{\lambda_2 \in \Lambda_2} \left|\hat f(\lambda_2)\right|^2
\geq p \sum_{\lambda\in\Lambda+\beta} \left|\tilde f_+ (\lambda)\right|^2 \tag 2.19
$$
for some positive constant $p$ (depending on the respective frame constants for the two g.s.p.'s). It will be convenient to write it in the form $p=2\alpha$, but this $\alpha$ will not be the one from \rom{(2.8)} in the theorem.

Returning to the summation \rom{(2.16)}, and still recalling the decomposition $f=f_+ + f_-$ from \rom{(2.11)}, we get
$$
\sum_{\lambda_1\in \Lambda_1} \left| \hat f(\lambda_1)\right|^2
= 2\alpha \sum_{\lambda\in\Lambda+\beta} |x_\lambda|^2
+\sum_{\lambda\in\Lambda+\beta} \left| (1+\overline{\langle\lambda,a \rangle}\,)x_\lambda
+(1-\overline{\langle \lambda,a\rangle}) y_\lambda \right|^2
$$
where the following shorthand notation was introduced:
$$
x_\lambda := \tilde f_+ (\lambda)\qquad\text{and}\qquad y_\lambda := \tilde f_-(\lambda)
\qquad\text{for }\lambda\in\Lambda+ \beta. \tag 2.20
$$

But this is the quadratic form which for each $\lambda \in \Lambda +\beta$ is given by the 2 by 2 matrix
$$
2\cdot \pmatrix \alpha + 1 +\cos(2\pi \lambda \cdot a)
& i\sin (2\pi \lambda \cdot a)\\
-i \sin (2\pi \lambda \cdot a)
& 1-\cos (2\pi \lambda \cdot a)\endpmatrix\tag 2.21
$$
and which by simple algebra is unitarily equivalent to the diagonal one $\pmatrix r_+ & 0\\ 0 &r_- \endpmatrix$ where the roots $r_\pm$ are given by
$$
r_\pm = 2+\alpha \pm \sqrt{(2+\alpha)^2 -4\alpha (1-\cos(2\pi \lambda \cdot a))}
$$
and this is formula \rom{(2.7)} where points in $\Lambda +\beta$ are written in the form $\lambda +\beta$ $(\lambda \in \Lambda)$.

If \rom{(i)--(ii)} both hold, it is clear that $r_\pm(\lambda)>0$ for all $\lambda \in \Lambda$, see formula \rom{(2.7)} in the theorem, and (on account of (ii)) as $\lambda$ varies in $\Lambda$, the set $\{r_\pm(\lambda) : \lambda \in \Lambda\}$ will then in fact be {\it finite}. Writing out the quadratic form \rom{(2.16)--(2.17)}, and using
$$
\langle f, F_1^*F_1 f\rangle = \sum_{\lambda\in\Lambda_1} |\hat f(\lambda_1)|^2
$$
for the sum in \rom{(2.16)} above, we conclude that there is an orthogonal decomposition of $\cl^2(\Omega_1)$ with corresponding projections $\{ E_\pm (\lambda)\}_{\lambda\in \Lambda}$ such that
$$
F_1^* F_1 = \sum_{\lambda\in\Lambda} r_+(\lambda)E_+(\lambda)
   + \sum_{\lambda\in \Lambda}r_- (\lambda)E_-(\lambda). \tag 2.22
$$
Notice that each of the non-zero projections (in the rational case)
$$
\sum_\lambda \{E_\pm (\lambda):\lambda \in \Lambda,\ r_\pm(\lambda) = \xi_\pm\}
$$
is infinite-dimensional. Moreover, the rationality condition (ii) guarantees that the numbers $\{r_- (\lambda) : \lambda \in \Lambda\}$ do not accumulate at $0$. This means that we get the lower estimate in condition \rom{(1.5)} for $(\Omega_1, \Lambda_1)$ being a generalized spectral pair. Conversely, if $0$ is in the closure of $\{r_-(\lambda) : \lambda \in \Lambda\}$ for some (fixed) $a$ and $\beta$, then this means that the lower bound estimate \rom{(1.5)} for $(\Omega_1, \Lambda_1)$ fails; and, by \rom{(2.7)} and \rom{(2.21)}, this means that one of the conditions (i), or (ii), must be violated.\qed\enddemo

We also get from the proof the following:

\proclaim{Corollary 2.1.3} Let the setup be as in the theorem, but assume instead that $(\Omega,\Lambda)$ and $(\Omega_2,\Lambda_2)$ are only g.s.p.'s. Then we still get that the g.s.p.\ property for $(\Omega_1,\Lambda_1)$ from \rom{(2.1)} and \rom{(2.6)} is equivalent to conditions \rom{(i)--(ii)}, but we will then \rom{not} have the exact formula for the frame operator $F_1^*F_1$ in \rom{(2.3)}, but instead of \rom{(2.22)} only the estimate
$$
F_1^*F_1 \geq \sum_{\lambda\in\Lambda} (r_+(\lambda)E_+(\lambda)
+r_-(\lambda)E_-(\lambda))\tag 2.23
$$
and the number $\alpha$ in \rom{(2.19)} and \rom{(2.7)} will depend on the given frame estimates for the two g.s.p.'s $(\Omega,\Lambda)$ and $(\Omega_2,\Lambda_2)$.\endproclaim

\demo{Proof} We have already read off the lower bounds from the proof above. The upper bound for the new pair $(\Omega_1, \Lambda_1)$ in \rom{(2.2)} and \rom{(2.6)} may be gotten from an application of \cite{JoPe7; Theorem 4.2}. \qed\enddemo

In summary, the difference between the respective cases, spectral pair and g.s.p., amounts to the difference between the equality in \rom{(2.22)} and the estimate in \rom{(2.23)}.

\remark{Remark 2.1.4} In Theorem 2.1.2 we construct new generalized spectral pairs \linebreak$(\Omega_1,\Lambda_1)$ from two given spectral pairs $(\Omega,\Lambda)$ and $(\Omega_2,\Lambda_2)$. We want to choose $\Lambda$ and a finite subset $B\subset \br^d$ such that $0\in B$, the non-overlapping property \rom{(1.12)} (generalizing \rom{(2.1)}) holds for $(\Omega,B)$, and $B\subset \Lambda_2^\circ$. Moreover $\Omega_2$ is extending $\Omega$, i.e., we assume the inclusion $\Omega \subset \Omega_2$. The theorem is about $B=\{0,a\}$ where $a$ is a single vector in $\br^d$, and here we point out the modifications which are entailed by the generalization to when $B$ contains more than two points, i.e., $|B| >2$. The above proof shows that when the frame operator
$$
F_1 : \cl^2(\Omega +B) \rightarrow \ell^2(\Lambda_1) \tag 2.24
$$
(where $\Lambda_1=\Lambda \cup \Lambda _2$ and $\Lambda$ is some spectrum for $\Omega$) is introduced, we then still have the formula \rom{(2.22)} above for $F_1^*F_1$ and the projections $E_\pm(\lambda)$ are indexed by $\Lambda$, and they are projections in $\ell_B^2$. However the eigenvalues $\{r_\pm (\lambda) : \lambda\in \Lambda\}$ will be given by a different formula (see details later). But it follows that we cannot have the spectrum of $F_1^* F_1$ bounded away from $0$. Specifically \rom{(2.22)} amounts to the decomposition
$$
F_1^*F_1 =\bigoplus_{\lambda\in\Lambda} Q(\lambda) \tag 2.25
$$
where each $Q(\lambda): \ell_B^2 \rightarrow \ell_B^2$ is a selfadjoint operator of rank 2. Since $\dim \ell_B^2 = |B|>2$, it follows that $Q(\lambda)$ cannot be strictly positive, and therefore, by \rom{(2.21)}, that $0$ must be in the spectrum of $F_1^*F_1$. This means, by \rom{(2.7)}, that we cannot have a positive lower bound in the frame estimate \rom{(1.5)} for $F_1$.

The argument for why each $Q(\lambda)$ must be of $\text{rank} \leq 2$ for each $\lambda\in\Lambda$ is based on a modification of the quadratic form \rom{(2.25)} in the proof of Theorem 2.1.2. When $B$ is given, introduce the function on $\br^d$:
$$
\varphi_B(\xi) := \sum_{b\in B} \langle b,\xi \rangle, \tag 2.26
$$
and note that $Q(\lambda)$ may be represented as a positive definite function on $\br^d \times \br^d$, or rather on $B^* \times B^*$ where $B^*$ is a subset of $\br^d$ with $|B^*| = |B|$, the positive definite function being
$$
Q_\lambda(\xi,\xi') 
= \overline{\varphi_B(\xi)}\,\varphi_B(\xi') \frac{m_d(\Omega_2)}{m_d(\Omega)}
+\overline{\varphi_B(\xi - \lambda)}\varphi_B (\xi'-\lambda) \tag 2.27
$$
for all $(\xi,\xi')\in B^* \times B^*$.\endremark

\subhead 2.2 Translation and scaling \endsubhead
In Section 1 we outlined the construction of new generalized spectral pairs (g.s.p.) on the form $(\Omega +B,\Lambda +L)$ where $(\Omega,\Lambda)$ is a given g.s.p.\ in $\br^d$, and $B$ and $L$ are finite subsets of $\br^d$ such that we have non-overlap of distinct translates $\Omega +b$ as $b$ varies over $B$, and similarly for $\{ \Lambda +\ell : \ell \in L\}$. In Section 2.1 we considered a different candidate for a possible spectrum of $\Omega + B$, but in this section we will give conditions for getting a spectrum of the form $\Lambda +L := \bigcup_{\ell\in L}(\Lambda +\ell)$, again with non-overlap for the distinct $\ell$-translates.

\proclaim{Theorem 2.2.1} Let $(\Omega,\Lambda)$ be a generalized spectral pair (g.s.p.) in $\br^d$ with frame constants $k$, $K$ as defined in \rom{(1.5)}, and let $B$, $L$ be two finite subsets in $\br^d$ satisfying the respective non-overlapping conditions and suppose
$$
B \subset \Lambda^\circ.\tag 2.28
$$
Let $\cm := (\langle b,\ell \rangle)$ be the corresponding matrix where $\langle b,\ell \rangle := e^{i 2\pi b\cdot \ell}$ are the entries, rows indexed by $L$ and columns by $B$. Suppose $\cm^*\cm$ is bounded below, or equivalently that it is invertible in $\ell_B^2$. With constants $p$, $P$ such that
$$
p\| v\| _{\ell_B^2} \leq \langle v, \cm^*\cm v \rangle_{\ell_B^2}
\leq P\|v\|_{\ell_B^2} \qquad\text{for }\forall v\in \ell_B^2,\tag 2.29
$$
we have $p>0$, and $(\Omega +B,\Lambda + L)$ is again a g.s.p.\  with frame constants multiplying, i.e., the respective lower and upper frame constants, $k_1$ respectively, $K_1$ are given by
$$
k_1 = kp \qquad\text{and}\qquad K_1=KP. \tag 2.30
$$\endproclaim

\demo{Proof} Consider four constants $A_1$, $A_2$, $k$, and $K$ such that $0<k\leq K$. We will be interested in double estimates on the form
$$
kA_2 \leq A_1 \leq KA_2, \tag 2.31
$$
and we will summarize this in the notation
$$
A_1 \cong (k,K)A_2. \tag 2.32
$$
From the assumptions on the g.s.p.\ $(\Omega,\Lambda)$ we get \rom{(2.32)} satisfied for
$$
A_1 =\sum_{\lambda \in \Lambda}|\tilde f(\lambda)|^2,
$$
and
$$
A_2 = \int_\Omega |f(x)|^2 \, dx = \|f\|_{\clo}^2. 
$$
The object is to establish that if $A_1 ' := \sum_{\Lambda + L} |\hat f(\lambda)|^2$,
$$
A_2' := \int_{\Omega +B} |f(x)|^2\,dx =\|f\|_{\cl^2(\Omega +B)}^2,
$$
and $f\in \cl^2(\Omega +B)$, then
$$
A_1' \cong (k_1,K_1)A_2' \tag 2.33
$$
where the new constants are given in \rom{(2.30)}. When the Fourier transform is on $\Omega +B$ it is denoted $\hat f$, and when it is calculated on $\Omega$ for the restriction, it is denoted $\tilde f$. Introducing the operator
$$
T_\ell f(x) = \sum_{b\in B} \overline{\langle\ell,b\rangle}\,f(x+b)\tag 2.34
$$
from $\cl^2(\Omega + B)$ to $\clo$ we then get
$$
\hat f (\lambda + \ell)=\widetilde{T_\ell f}\,(\lambda + \ell),
$$
and
$$
\sum_{\lambda\in\Lambda}\left| \widetilde{T_\ell f}\,(\lambda +\ell)\right|^2
\cong (k,K)\times \|T_\ell f\|_{\clo}^2\tag 2.35
$$
for all $f\in \cl^2(\Omega +B)$ and all $\ell \in L$. Recall that each $\ell + \Lambda$ serves as a spectrum for $\Omega$ and the frame bound is the same.

For the second summation we get, using assumption \rom{(2.29)},
$$\align
\sum_\ell \int _\Omega |T_\ell f(x)|^2 \,dx
&=\sum_b \sum_{b'} \sum_\ell \langle \ell, b-b' \rangle
   \int_\Omega \overline{f(x+b)}\,f(x+b')\,dx\\
&\cong (p,P)\int_\Omega \sum_b |f(x+b)|^2 \, dx\\
&=(p,P) \int_{\Omega +B} |f|^2\,dx\\
&=(p,P) \| f\|^2 _{\cl^2(\Omega + B)},
\endalign$$
and the result now follows by substitution.\qed\enddemo

\remark{Remark} Systems of operators like \rom{(2.34)} abound in wave-packet analysis, and in wavelet theory, see e.g., \cite{JoPe3--6} and \cite{Mey}. The corresponding adjoint operators, $T_\ell^* : \clo \rightarrow \cl^2 (\Omega +B)$, are given by
$$
T_\ell^* f(x) = \sum _{b\in B} \chi_{\Omega+b} (x)
\langle \ell,b\rangle f(x-b)
\qquad\text{for $\forall f\in \clo$},
$$
and it follows that
$$
T_{\ell'} T_\ell^* = \varphi_B(\ell ' - \ell) I_{\clo}\qquad\text{for $\forall \ell$, $\ell' \in L$},
$$
where $\varphi_B$ is given in \rom{(2.26)}, i.e.,
$$
\varphi_B(\xi):= \sum_{b\in B} \langle b,\xi \rangle
= \sum_{b\in B} e^{i2\pi b\cdot \xi}\qquad\text{for $\forall \xi \in \br^d$}.
$$
In the special case where
$\varphi_B (\ell' - \ell) =|B|\delta_{\ell,\ell'}$, we then get an associated {\it orthogonal\/} decomposition of $\cl^2(\Omega+B)$, corresponding to the disjoint decomposition, $\bigcup_{b\in B}(\Omega+b)$, into translates of a fixed domain $\Omega$ (null-overlap of distinct translates is allowed). A main point in the present paper is to show that wave-packet analysis carries over to the almost orthogonal case.
\endremark

The next result is about {\it stability\/} of the exactness property (see Section 1) in the step
$$
(\Omega,\Lambda)\mapsto (\Omega +B,\Lambda +L) \tag 2.36
$$
in the category of g.s.p. Recall {\it exactness\/} for a given g.s.p.\ $(\Omega,\Lambda)$ is the property that the associated operator $F: \clo \rightarrow \ell^2(\Lambda)$ maps onto $\ell^2(\Lambda)$. So it is equivalent to the adjoint
$$
F^* : \ell^2(\Lambda) \rightarrow \clo
$$
having zero kernel, see \rom{(1.7)} above.

\proclaim{Corollary 2.2.2} If the initial g.s.p.\ $(\Omega,\Lambda)$ in the construction \rom{(2.36)} from Theorem 2.2.1 is exact, and if in addition to $B\subset \Lambda^\circ$, it is assumed that the matrix $\cm=(\langle b,\ell \rangle)$ is invertible, i.e., $|B|=|L|$ and $\cm$ invertible as a square matrix, then it follows that the new g.s.p.\ $(\Omega +B,\Lambda +L)$ is also exact.\endproclaim

\demo{Proof (sketch)} Let $(\xi)\in \ell^2 (\Lambda +L)$ and assume that $F_1^* \xi=0$ in $\cl^2(\Omega +B)$ where $F_1$ is the frame operator for the new g.s.p.\ in the Theorem 2.2.1 construction. Then the double sum
$$
\sum_{\ell\in L} \sum_{\lambda\in \Lambda} \xi_{\ell + \lambda} 
\langle \ell +\lambda, x+b \rangle
$$
must vanish a.e.\ $x\in \Omega$, and for all $b\in B$. Since $B\subset \Lambda^\circ$, we get
$$
\sum_{\ell \in L} \langle \ell,b \rangle \sum_{\lambda\in\Lambda} \xi_{\ell + \lambda}
\langle \ell + \lambda, x \rangle =0.
$$
Invertibility of $\cm$ yields
$$
\sum_{\lambda\in \Lambda} \xi_{\ell + \lambda} 
\langle \ell + \lambda, x \rangle =0
$$
for all $\ell \in L$ and a.e.\ $x\in \Omega$. Using exactness of $(\Omega,\Lambda)$ we conclude that $\xi_{\ell + \lambda}=0$ for all $\ell\in L$ and all $\lambda \in \Lambda$.\qed\enddemo

The final result is about the multiplicative property. Let $R$ be a $d$ by $d$ matrix over the reals and assume that all its eigenvalues satisfy $|\lambda_i |> 1$. The transpose is denoted $R^*$. For the lemma, we just need $R$ to be invertible, but the stronger property will be needed in the applications in the next section.

\proclaim{Lemma 2.2.3} Suppose $(\Omega,\Lambda)$ is a g.s.p.\ in $\br^d$ and $R$ is as specified. Then $(R^{-1}\Omega, R^*\Lambda)$ is also a g.s.p.\ and the frame constants $(k,K)$ are the same for the two g.s.p.'s.\endproclaim

\demo{Proof} This simple argument is based on the transform formula for the Lebesgue measure in $\br^d$ and we leave details to the reader.\qed\enddemo

\subhead 2.3 Reversing the construction \endsubhead
We have seen how to get new generalized spectral pairs from old ones by the step $\Omega \mapsto \Omega +B$ where $\Omega$ has a spectrum in the generalized sense and when $B$ is a finite set such that the translates $\{\Omega +b: b\in B\}$ do not overlap. In this section we show that the process can be reversed but this reversion is surprisingly subtle, and it is relevant in tiling theory, see \cite{JoPe3--4} and especially \cite{LaWa}. We show that if $B$ has the non-overlap property and if there is a set $\Lambda$ such that $(\Omega +B, \Lambda)$ is a generalized spectral pair, then there is an explicit construction of a second set $\Lambda_\Omega$ such that $(\Omega,\Lambda_\Omega)$ is then also a generalized spectral pair. But as pointed out in \cite{LaWa} even if $d=1$ and $\Omega + B$ is a $\bz$ fundamental domain, i.e., if $(\Omega +B,\bz)$ is a {\it spectral pair\/} (in the strict sense) in $\br$, then typically $(\Omega,\Lambda_\Omega)$ will still only be a {\it generalized\/} spectral pair.

\proclaim{Theorem 2.3.1} 
\roster
\item"{(i)}" Let $\Omega$, $\Lambda$, $B$ be subsets of $\br^d$, with $B$ a finite set, and $0\in \Lambda$. Suppose $B+\Omega =\bigcup_{b\in B} (b+\Omega)$ is non-overlapping, and that $(B+\Omega, \Lambda)$ is a generalized spectral pair. Assume $\langle e_0,e_\lambda \rangle_{\Omega +B} =0$ for all $\lambda \in \Lambda\backslash \{0\}$. Let
$$
\Lambda_\Omega := \{\lambda\in \Lambda : \hat \chi_\Omega (\lambda)=0\} \cup \{0\}.
$$
Then $(\Omega, \Lambda_\Omega)$ is a generalized spectral pair.
\item"{(ii)}" If $(B+\Omega, \Lambda)$ has frame constants $(k,K)$, then $(\Omega,\Lambda_\Omega)$ has frame constants $(k/|B|,K)$.\endroster\endproclaim

\demo{Proof} For functions $g$ on $\br^d$, let $Z(g)=\{ \lambda\in\Lambda : g(\lambda)=0\}$, then $\Lambda_\Omega =Z(\hat \chi_\Omega )\cup (0)$. Let $\lambda \in \br^d$, then 
$$
\hat\chi_{B+\Omega}(\lambda) 
  = \int_{B+\Omega} \overline{\langle t,\lambda \rangle}\, dt
  = \sum_{b\in B} \overline{\langle b,\lambda \rangle} \int_\Omega \overline{\langle x,\lambda \rangle}\, dx
 = \varphi_B (\lambda)\hat\chi_\Omega(\lambda).
$$
Hence, $\Lambda\backslash (0) =Z(\hat\chi_{B+\Omega}) = Z(\varphi_B)\cup Z(\hat\chi_\Omega)$; in particular $\varphi_B(\lambda) =0$ for all $\lambda\in \Lambda \backslash \Lambda_\Omega$. Let $f\in \clo$, and define $g\in \cl^2 (B+\Omega)$ by $g(b+x)=f(x)$ for $b\in B$, $x\in \Omega$ (i.e., $g=T_0 f$ in the notation used above).
$$\align
k|B| \,\|f\|_\Omega^2 
  &= k\| g\|_{B+\Omega}^2 \leq \sum_{\lambda\in\Lambda} |\hat g(\lambda)|^2\\
  &=\sum_{\lambda\in\Lambda} \left| \sum_{b\in B}\int_{\Omega} 
      \overline{\langle b,\lambda\rangle}\, \overline{\langle x,\lambda\rangle}\,
      f(x)\,dx \right|^2\\
  &=\sum_{\lambda\in\Lambda} \left| \varphi_B (\lambda)\tilde f(\lambda)\right|^2\\
  &= \sum_{\lambda\in\Lambda_\Omega} \left|\varphi_B (\lambda)\tilde f(\lambda)\right|^2\\
  &\leq |B|^2 \sum_{\lambda\in\Lambda_\Omega} |\tilde f(\lambda)|^2.
\endalign$$
This gives the lower bound. The upper bound follows from $\Lambda_\Omega \subset \Lambda$.\qed\enddemo

There is a much stronger version of the assumptions in Theorem 2.3.1 which stills gives an interesting application to the classification of {\it spectral pairs\/} in $\br^d$. First note that a subset $\cd\subset\br^d$ is a {\it fundamental domain\/} if it is measurable and there is a lattice $\Gamma$ in $\br^d$ such that $\cd +\Gamma =\br^d$ with non-overlapping translates $\{\cd + \gamma : \gamma\in \Gamma\}$. We say that $(\cd,\Gamma)$ is a {\it lattice tiling}. In \cite{JoPe3} we considered subsets in $\br^d$ on the form $\cd +B$ with $B$ finite and non-overlap for the translates $\{\cd+b : b\in B\}$. A main question is to decide the structure of spectral pairs on the form $(\cd+B,\Lambda)$  when the sets $\cd$ and $B$ are specified as stated. (Note that if $B\subset \Gamma$, then the non-overlap property will be automatic.) If then $0\in \Lambda$ we clearly have the assumption,
$$
\int_{\cd+B} e_\lambda\, dx = 0\qquad \text{for } \forall \lambda\in \Lambda\backslash \{0\},
\tag 2.37
$$
satisfied. But if, as in Theorem 2.3.1, $(\cd +B,\Lambda)$ is only assumed to be a g.s.p., then this may not be the case. In Example 2.1.1 above, \rom{(2.37)} is satisfied in fact iff $a=1$ which is really the degenerate case. The following result serves as a converse to \cite{JoPe3; Theorem 6.1} and classifies an important subclass of spectral pairs in~$\br^d$.

\proclaim{Theorem 2.3.2} Let $\cd$ be a fundamental domain in $\br^d$, and let $\Lambda$, $B$ be subsets in $\br^d$ such that $B$ is finite and satisfies the non-overlap property relative to $\cd$, $0\in \Lambda$, and $(\cd + B,\Lambda)$ is a spectral pair in $\br^d$. Suppose that there is some lattice $\Gamma$ such that $(\cd,\Gamma)$ is a lattice tiling, and that we have the following inclusion
$$
\{\lambda\in\Lambda : \hat\chi_\cd (\lambda)=0\} \subset \Gamma^\circ.
$$
Then it follows that there is a second finite subset $L\subset \br^d$ of same cardinality as $B$ such that $\Lambda = L+\Gamma^\circ$ where $\Gamma^\circ$ denotes the lattice which is dual to $\Gamma$.\endproclaim

\demo{Proof} From Theorem 2.3.1 we conclude that $(\cd,\Lambda_\cd)$ is a g.s.p.\ where
$$
\Lambda_\cd =\{ \lambda\in \Lambda : \hat\chi_\cd (\lambda) =0\} \cup \{ 0 \}. \tag 2.38
$$
But since $(\cd,\Gamma)$ is a lattice tiling, and the pair $(\cd+B,\Lambda)$ is a spectral pair, we conclude for the set $\Lambda_\cd$ in \rom{(2.38)} that $\Lambda_\cd \subset \Gamma^\circ \cap \Lambda$. So $\Lambda_\cd \subset \Gamma^\circ$ and $\{e_s|_\cd : s\in \Gamma^\circ\}$ is an orthogonal basis for $\cl^2(\cd)$. Since $\{e_s|_\cd : s\in \Lambda_\cd\}$ is total in $\cl^2(\cd)$, (recall that the pair $(\cd,\Lambda_\cd)$ was a g.s.p.\ by Theorem 2.3.1), we conclude that
$$
\Lambda_\cd = \Gamma^\circ \cap \Lambda = \Gamma^\circ,
$$
and so $\Gamma^\circ \subset \Lambda$. If $\ell_1 \in \Lambda \backslash \Gamma^\circ$ we may apply the same argument and Theorem 2.3.1 to the spectral pair $(\cd +B, \Lambda - \ell_1)$ and conclude that $\Gamma^\circ \subset \Lambda - \ell_1$, or equivalently $\{0,\ell_1\} + \Gamma^\circ \subset \Lambda$. If $\ell_2 \in \Lambda \backslash (\{0,\ell_1\}+\Gamma^\circ)$, then the same argument would give $\{0,\ell_1,\ell_2\} + \Gamma^\circ \subset \Lambda$. But an application of \cite{JoPe3; Theorem 6.1} and \cite{Lan; Theorem 4} shows that the inductive process must stop, and then the corresponding matrix $\left( |B|^{-1/2} e^{i 2\pi \ell_j \cdot b_k} \right)_{j,k}$ will be unitary. So in particular the set $L:=\{0,\ell_j\}$ will have the same cardinality as $B$, and $L+\Gamma^\circ =\Lambda$.\qed\enddemo

Further details and applications of this result to the classification problem for spectral pairs in $\br^d$ will be given in a subsequent paper by the coauthors.

\subhead 2.4 Spectral pairs in one dimension \endsubhead
This section serves as an appendix to Section 2.3 above, and shows that, when $d=1$, and when the fundamental domain $\cd$ from Theorem 2.3.2 is further assumed to be the unit-interval $I=(0,1)$, then we have all the extra assumptions in Theorem 2.3.2 automatically satisfied. As a corollary, we therefore get a classification of all spectral pairs (for $d=1$) on the form $(\Omega,\Lambda)$ where $\Omega$ is a finite (non-overlapping) union of translates of $I$. Specifics are worked out in detail for the convenience of the reader.

The next result, as the one before it, has a {\it converse}, which in fact follows from our earlier paper \cite{JoPe3; Sections 6--8}. Putting the results from the two papers together, we then get a {\it complete classification\/} of all spectral pairs $(\Omega,\Lambda)$ in one dimension where the set $\Omega$ is built from translates of an interval; specifically, $\Omega = I+B = \bigcup_{b\in B} (I+b)$, with non-overlap, where $I=(0,1)$, and $B\subset \br$ is a finite subset.

\proclaim{Corollary 2.4.1} Let $I=(0,1)$, $B\subseteq \br$ be finite, $B+I$ be non-overlapping, and let $\Lambda\subseteq \br$. Assume that $(B+I,\Lambda)$ is a spectral pair. Then $\Lambda =L +\bz$ for some finite set $L$, $|L|=|B|$, $L+\bz=\bigcup_{\ell\in L} \ell+\bz$ is non-overlapping, and $b-b' \in \bz$ for all $b$, $b'\in B$. Moreover, the matrix $(|B|^{-1/2} \langle b,\ell \rangle )_{b\in B, \ell \in L}$ must be unitary.
\endproclaim

\demo{Proof} Let $\ell\in \Lambda$, and define 
$$
\Lambda_I := \{\lambda\in \Lambda-\ell : \hat\chi_I (\lambda)=0\} \cup \{0\},
$$
where
$$
\hat\chi_I (t) = \int_I e_t (x)\,dx
=\cases 1,&\text{if $t=0$}\\ \frac{e^{i2\pi t}-1}{i2\pi t}, &\text{if $t\neq 0$.}\endcases
$$
Hence $\hat\chi_I(t) = 0$ iff $t\in\bz\backslash (0)$, and $\Lambda_I\subseteq \bz$. By Theorem 2.3.1, $(I,\Lambda_I)$ is a generalized spectral pair; in particular, $\{e_\lambda :\lambda \in \Lambda_I \}$ is total in $\cl^2(I)$. But $\{ e_\lambda :\lambda\in \bz\}$ is an ONB for $\cl^2(I)$, and $\Lambda_I \subset \bz$; thus $\Lambda_I = \bz$. It follows that $\ell + \bz = \ell + \Lambda_I \subset \Lambda$. Inductively, suppose $\{\ell_1,\ldots, \ell_p\} + \bz \subset \Lambda$, and let $\ell\in \Lambda \backslash(\{\ell_1,\ldots,\ell_p\} +\bz)$. Arguing as above, we get $\{\ell, \ell_1,\ldots, \ell_p \} + \bz \subset \Lambda$. By \cite{Lan, Theorem 4} we have 
$$
d^+(\Lambda) = \limsup_{r\rightarrow\infty}
 \left(r^{-1}\max_{a\in \br} |\Lambda\cap B_a (r)|\right)
\leq m_1 (\Omega)
$$
where $B_a(r) = [a,a+r]$. Hence the construction above must terminate, and $\Lambda = L+\bz$ for some finite set $L\subset \Lambda$. By construction, $\ell - \ell' \notin \bz$ for all $\ell$, $\ell' \in L$ with $\ell\neq\ell '$. Let
$$
\varphi_B (t) = \sum_{b\in B} \langle b,t \rangle,\qquad t\in \br.
$$
Then 
$$
\varphi_B (\ell - \ell ' + z-z') \cdot \hat\chi_I (\ell -\ell ' +z -z')
= \int_{B+I} \overline{e_{\ell'+z'}} \, e_{\ell + z} =0
$$
if $\ell ' + z' \neq \ell + z$. It follows that
$$
\varphi_B (\ell - \ell' + z-z')=0,\qquad\text{for all $z$, $z'\in\bz$, $\ell$, $\ell'\in L$}
$$
with $\ell\neq \ell '$. For $z\in \bz$, let $\cm_z$ be the $|B| \times |L|$ matrix
$$
\cm_z =(\langle b,\ell +z\rangle )_{b\in B,\ell\in L}
$$
($b\in B$ labelling the rows, and $\ell \in L$ labelling the columns). Then $\cm_z$ has orthogonal columns; hence $|B|\geq |L|$. By totality of $\{e_{\ell+z} : \ell\in L,\ z\in\bz\}$ in $\cl^2(B+I)$, we have $|L| \leq |B|$. Thus $|B|=|L|$, and $\cm_z$ is a Hadamard matrix for all $z\in \bz$.

For the rest of the proof, we will without loss of generality assume that $0\in B$ and $0\in L$. Define vectors $v_n$, $w_\ell \in \ell_B^2$ by $v_n = (\langle b,n \rangle)_{b\in B}$, $w_\ell = (\langle b,\ell \rangle )_{b\in B}$ for $n\in \bz$ and $\ell \in L$. Then
$$
v_n \cdot w_\ell = \sum_b \langle b, \ell -n \rangle =\varphi_B (\ell -n).
$$
Hence $v_n \cdot w_\ell =0$ for $\ell \in L\backslash (0)$, $n\in \bz$. The vectors in $\{w_\ell : \ell\in L\backslash (0)\}$ are mutually orthogonal, and hence span an $(|L|-1)$-dimensional subspace in $\ell_B^2$. It follows that there are scalars $c_n \in\bc$ such that
$$
v_n = c_n v_0=c_n w_0
$$
for all $n\in\bz$. (Note $v_0 = (1,1,\ldots,1)$.)

So $\langle b,n \rangle$ is independent of $b\in B$. Therefore $\langle b,1 \rangle = \langle 0,1 \rangle =1$, and it follows that $b\in B$; i.e., $B\subseteq \bz$.\qed\enddemo

\head 3. The Spectrum of the Limiting Fractal \endhead
\subhead 3.1 Iteration \endsubhead
In this section we consider results and examples which realize limit properties for the approximating discrete transforms
$$
\cf_n f=\widehat{f\,d\mu_n} = \int \bar e_\lambda f\,d\mu_n
$$
which enter into our iteration theory as outlined in Section 1 above. For a given generalized spectral pair $(\Omega,\Lambda)$ in $\br^d$ we shall need Lemma 1.1 about the discrete set $K := \Lambda^\circ$; and we shall also need the simple observation that every translate $(\Omega +a,\Lambda +\beta)$ is again a g.s.p.\ with the same frame constants as the original one $(\Omega,\Lambda)$.

Let $(\Omega,\Lambda)$ be a generalized spectral pair, let $B$ and $L$ be finite subsets of $\br^d$, let $R$ be an invertible $d\times d$ matrix with real entries. Define affine maps of $\br^d$, $\sigma_b$ and $\tau_\ell$, for $b\in B$ and $\ell\in L$ by
$$
\sigma_b x=R^{-1} (x+b)\qquad\text{and}\qquad\tau_\ell s=R^*(s+\ell) \tag 3.1
$$
for all $x$, $s\in\br^d$. Inductively define sequences $(\Omega_n)$, $(\Lambda_n)$, $(n=0,1,\ldots)$ of sets by 
$$\align
\Omega_0& = \Omega,\qquad\Lambda_0=\Lambda, \tag 3.2\\
\Omega_{n+1} &= \sigma(\Omega_n) =\bigcup_{b\in B} \sigma_b(\Omega_n),\tag 3.3\\
\intertext{and}
\Lambda_{n+1} &=\tau(\Lambda_n) =\bigcup_{\ell\in L}\tau_\ell (\Lambda_n).\tag 3.4
\endalign$$

The following five conditions will be standing assumptions for our duality \linebreak iteration-algorithm:
\roster
\item"{(i)}" the matrix $\cm=(\langle \ell,b\rangle)_{\ell\in L,b\in B}$, where $\langle \xi,t \rangle = e_\xi(t)=e^{i2\pi\xi \cdot t}$, has rank $|B|$, 
\item"{(ii)}" $R(K)\subset K$,
\item"{(iii)}" $B\subset K$,
\item"{(iv)}" $L\subset R^{*-1} (K^\circ)$, and
\item"{(v)}" the sum $L +K^\circ$ is direct; i.e., $(\ell+K^\circ) \cap (\ell'+K^\circ)=\emptyset$ for $\forall \ell \neq \ell'$ in $L$.\endroster
Here $K=\Lambda^\circ$ as before. Then $(\Omega_n,\Lambda_n)$, $n=0,1,\ldots$, is a sequence of generalized spectral pairs, as will be shown below. In addition to the results from Section 2 above we shall also need the following four lemmas:

\proclaim{Lemma 3.1.1} Let $\Omega$, $\Lambda$, and $R$ be given as described above, and assume that conditions \rom{(i)--(v)} hold. Let $(\Omega_n)$ and $(\Lambda_n)$ be the iterations. We then have: For each $n\in \{ 0,1,2,\ldots\}$
\roster
\item"{(a)}" $B\subset \Lambda_n^\circ$.
\item"{(b)}" $\sigma_b(\Omega_n) \cap \sigma_{b'}(\Omega_n)$ is a null-set whenever $b$, $b'\in B$ and $b\neq b'$.
\item"{(c)}" $\tau_\ell(\Lambda_n) \cap \tau_{\ell'}(\Lambda_n) = \emptyset$ whenever $\ell$, $\ell'\in L$ and $\ell\neq\ell '$.\endroster\endproclaim

\demo{Proof} 
\roster\item"{(a)}" Recall $\Lambda_n =\{R^{*n}\lambda + \sum_{k=1}^n R^{*k} \ell_k : \lambda \in \Lambda,\ \ell_k\in L\}$. Using \rom{(ii)}, $K=\Lambda^\circ$, \rom{(iv)}, and \rom{(iii)}, we get $R^*K^\circ \subset K^\circ$, $\Lambda\subset \Lambda^{\circ\circ}=K^\circ$, $R^*L\subset K^\circ$, and $K^\circ \subset B^\circ$. Hence, $\Lambda_n\subset K^\circ \subset B^\circ$. Taking polars in turn yields, $B\subset B^{\circ\circ}\subset K \subset \Lambda_n^\circ$.

\item"{(c)}" From the proof of (a), we have $\Lambda_n \subset K^\circ$; hence the desired conclusion follows from the asumption \rom{(v)} that $(K^\circ + \ell) \cap (K^\circ + \ell ') = \emptyset$.

\item"{(b)}" Claim: If $\{e_\rho : \rho \in \cR\}$ is total in $\cl^2(S)$, then $S \cap (S+\eta)$ is a null-set for any $\eta\in \cR^\circ\backslash(0)$. Let $f(x)=1$ if $x\in S(\eta)=S \cap (S+\eta)$, and $f(x)=0$ if $x\in S\backslash S(\eta)$. Note that $S(\eta)\subseteq S$, $S(\eta)-\eta\subseteq S$, and $e_\rho(x)=e_\rho(x-\eta)$ for all $x\in S(\eta)$ and all $\rho\in \cR$. But $f(x)=1\neq 0 = f(x-\eta)$ when $x\in S(\eta)$, so the totality of $\{e_\rho: \rho \in \cR\}$ implies that $\eta=0$, or $S(\eta)$ is a null-set. This proves the claim. \endroster

{\it Completion of the proof of \rom{(b)}.} We will show, by induction on $n$, that $(\Omega_n,\Lambda_n)$ is a generalized spectral pair, and that $\sigma_b(\Omega_n) \cap \sigma_{b'}(\Omega_n)$ is a null-set. For $n=0$, this follows from the assumption that $(\Omega,\Lambda)$ is a generalized spectral pair, and the claim.

If $(\Omega_n,\Lambda_n)$ is a generalized spectral pair and the intersections $\sigma_b(\Omega_n) \cap \sigma_{b'}(\Omega_n)$ are null-sets, then it follows from part \rom{(a)}, part \rom{(c)}, and Theorem 2.2.1 and Lemma 2.2.3, that $(\Omega_{n+1}, \Lambda_{n+1})$ is a generalized spectral pair, but then the claim, together with the observation $\Lambda_{n+1} \subset K$ (from the proof of part \rom{(a)}), finishes the proof.\qed\enddemo

Let $\Delta$ denote a Borel set in $\br^d$. Let $m(\Delta)$ denote the Lebesgue measure of the set $\Delta$. Using Lemma 3.1.1(b), we then note that
$$
\Omega_n =\left\{ \sum_{k=1}^n R^{-1} b_k : b_k \in B\right\} + R^{-n}\Omega \tag 3.5
$$
is without overlap (redundancy), so we may define a probability measure $\mu_n$ on $\br^d$~by
$$
\mu_n(\Delta) = \frac{|\det R|^n}{|B|^n m_d(\Omega)} m_d(\Delta \cap \Omega_n) \tag 3.6
$$
for Lebesgue measurable sets $\Delta$.

\proclaim{Lemma 3.1.2} The measures $\mu_n$ defined above satisfy
$$
\mu_{n+1} = \frac{1}{|B|} \sum_{b\in B} \mu_n \circ \sigma_b^{-1} \tag 3.7
$$
for $n=0,1,\ldots$.\endproclaim

\demo{Proof} Using Lemma 3.1.1(b) we get
$$\align
\int f \, d\mu_{n+1} 
  &= \frac{|\det R|^{n+1}}{|B|^{n+1}m(\Omega)} \int_{\Omega_{n+1}} f(x)\,dx\\
  &=\frac{|\det R|^{n+1}}{|B|^{n+1} m(\Omega)}
        \sum_{b\in B} \int_{\sigma_b(\Omega_n)} f(x)\,dx\\
  &=\frac{|\det R|^n}{|B|^{n+1}m(\Omega)} 
        \sum_{b\in B}\int_{\Omega_n} f(\sigma_b x)\, dx\\
  &=\frac{1}{|B|} \sum_{b\in B} \int f\circ \sigma_b \,d\mu_n.\qed
\endalign$$
\enddemo

If all the eigenvalues $\lambda_i$ of $R$ have $|\lambda_i|>1$, then Hutchinson's theorem \cite{Hut} tells us that the sequence $(\mu_n)$ converges to some limit probability measure $\mu$ satisfying $\mu=\frac{1}{|B|} \sum_b \mu\circ\sigma_b^{-1}$, or equivalently, $\int f \, d\mu = \sum_b \frac{1}{|B|} \int f\circ \sigma_b \, d\mu$ for all Borel functions $f$; and we now turn to the analysis of $\mu$.

For future reference we first record the simple

\proclaim{Lemma 3.1.3} If the constants $a\leq A$ are such that
$$
a\|f\|_\Omega^2 \leq \sum_{\lambda\in\Lambda} |\hat f(\lambda)|^2
\leq A\| f\|_\Omega^2 \tag 3.8
$$
for all $f \in \cl^2(\Omega)$, then
$$
\frac{a}{m_d(\Omega)} \| f\|_{\mu_0}^2
 \leq \sum_{\lambda\in\Lambda_0} |\cf_0 f(\lambda)|^2 
 \leq \frac{A}{m_d (\Omega)} \| f\|_{\mu_0}^2 \tag 3.9
$$
for all $f\in \cl^2(\mu_0)$.\endproclaim

From now on  we will assume that a triple, $\Lambda$, $B$, $L\subseteq \br^d$, a sequence of probability measures $\mu_n$ on $\br^d$, and an invertible $d\times d$ matrix $R$ are given. We will suppose that the five conditions \rom{(i)--(v)} are satisfied, and that
$$
\mu_{n+1} = \frac{1}{|B|} \sum_{b\in B} \mu_n \circ \sigma_b^{-1} \tag 3.10
$$
for $n=0,1,\ldots$. Notice that Lemma 3.1.1(a),(c) remain valid, with the same proof.

If all eigenvalues of $R$ have modulus $>1$, the theorem in \cite{Hut} tells us that $\mu_n$ converges to a probability measure $\mu$ satisfying $\mu =\frac{1}{|B|} \sum_{b\in B} \mu \circ \sigma_b^{-1}$.

\proclaim{Lemma 3.1.4} For $f\in \cl^2 (\mu_n)$, let 
$$
\cf_n f(\lambda) =\int f(x) \overline{\langle x,\lambda \rangle}\, d\mu_n (x).
$$
Then for $f\in \cl^2(\mu_{n+1})$, $\ell\in L$, $\lambda\in \Lambda_n$, we have the transformation rule
$$
\cf_{n+1} f(\tau_\ell \lambda ) 
= \cf_n \left(\bar e_\ell \frac{1}{|B|} \sum_{b\in B} 
  \overline{\langle b,\ell \rangle}\, f\circ \sigma_b\right)(\lambda).
$$
Recall that $\langle x,\lambda\rangle = e_\lambda(x)=e^{i2\pi\lambda\cdot x}$.\endproclaim

\demo{Proof}
$$\align
(\cf_{n+1} f)(\tau_\ell \lambda)
  &=\int  f(x)\,\overline{\langle x,\tau_\ell \lambda \rangle}\,d\mu_{n+1} (x)\\
  &= \frac{1}{|B|} \sum_{b\in B} \int f(\sigma_b x)
           \overline{\langle \sigma_bx,\tau_\ell\lambda \rangle}\,d\mu_n(x)\\
  &=\frac{1}{|B|} \sum_{b\in B} \overline{\langle b,\ell \rangle} 
          \int \overline{\langle x,\ell \rangle}\, f(\sigma_bx)\,
            \overline{\langle x,\lambda \rangle}\, d\mu_n(x)
\endalign$$
as desired. We used that $\langle b,\lambda \rangle =1$; this follows from Lemma 3.1.1(a). \qed\enddemo

\proclaim{Theorem 3.1.5} Let $\Lambda$, $L$, $B$, $R$, $\mu_n$ be as defined above. Assume there are constants $k\leq K$ such that
$$
k\| f\|_{\mu_0}^2 \leq \sum_{\lambda\in\Lambda_0} |\cf_0 f(\lambda)|^2
\leq K\| f\|_{\mu_0}^2 \tag 3.11
$$
for all $f\in \cl^2(\mu_0)$. Let $\cm = (\langle b,\ell\rangle )_{b\in B,\ell \in L}$, and let $d\leq D$ be such that the $|B|$ by $|B|$ matrix $\cm^*\cm$ satisfies 
$$
d\| v\|_{\bc^{|B|}}^2 \leq \langle v, \cm^*\cm v\rangle_{\bc^{|B|}}
\leq D\|v\|_{\bc^{|B|}}^2
$$
for all $v\in \bc^{\|B|}$. Then
$$
k\left( \frac{d}{|B|}\right)^n \| f\| _{\mu_n}^2
\leq \sum_{\lambda\in\Lambda_n} |\cf_n f(\lambda)|^2
\leq K\left( \frac{D}{|B|}\right)^n \| f\| _{\mu_n}^2 \tag 3.12
$$
for all $f\in \cl^2(\mu_n)$ and each $n=0,1,\ldots$.
\endproclaim

Here, $\Lambda_0=\Lambda$ and $\Lambda_{n+1} =\tau(\Lambda_n)$, as before, and $\cf_n$ is as in Lemma 3.1.4. Notice, by Lemma 3.1.3, that if $(\Omega,\Lambda)$ is a generalized spectral pair and $\mu_0 = \frac{1}{m(\Omega)} m$ (meaning $\mu_0 (\Delta) = \frac{1}{m(\Omega)} m(\Delta \cap \Omega)$ for all Borel sets $\Delta$). Then the hypothesis of this theorem is satisfied.

\demo{Proof} The proof is by induction on $n$. Assume that
$$
k\left( \frac{d}{|B|}\right)^n \| f\|_{\mu_n}^2 
\leq \sum_{\lambda\in\Lambda_n} |\cf_n f(\lambda) |^2 
\leq K\left( \frac{D}{|B|}\right)^n \| f\|_{\mu_n}^2
$$
for all $f\in \cl^2 (\mu_n)$.

For $f\in\cl^2 (\mu_{n+1})$, let
$$
T_\ell f = \bar e_\ell \frac{1}{|B|} \sum_{b\in B} \overline{\langle b,\ell \rangle} \, f\circ \sigma_b,
$$
then 
$$\align
\sum_{\ell\in L} \| T_\ell f\|_{\mu_n}^2
  &= \frac{1}{|B|^2} \sum_{\ell\in L} 
\left\| \sum_{b\in B} \overline{\langle b,\ell \rangle} f\circ \sigma_b\right\|_{\mu_n}^2\\
  &= \frac{1}{|B|^2} \sum_{b\in B}\sum_{b'\in B} 
 \left\langle f\circ \sigma_b,\sum_\ell \langle b-b',\ell \rangle             f\circ\sigma_{b'}\right\rangle_{\mu_n}\\
  &=\frac{1}{|B|^2} \int \langle v(x), \cm^*\cm v(x) \rangle _{\bc^{|B|}} \, d\mu_n (x)
\endalign$$
where $v(x)$ is the vector in $\bc^{|B|}$ with entries $(f \circ \sigma_b(x))_{b\in B}$. Hence $$\align
\frac{d}{|B|^2} \int \| v(x)\|_{\bc^{|B|}}^2 \,d\mu_n (x)
 &\leq \sum_{\ell \in L} \| T_\ell f\|_{\mu_n}^2\\
 & \leq \frac{D}{|B|^2} \int\| v(x) \|_{\bc^{|B|}} ^2 \, d\mu_n (x)
\endalign$$
it follows that
$$
\frac{d}{|B|} \| f\|_{\mu_{n+1}} ^2 
\leq \sum_{\ell\in L} \| T_\ell f \|_{\mu_n}^2
 \leq \frac{D}{|B|} \| f\|_{\mu_{n+1}}^2. \tag *
$$
Now
$$\align
\sum_{\xi\in\Lambda_{n+1}} |\cf_{n+1} f(\xi)|^2
&= \sum_{\ell\in L} \sum_{\lambda\in \Lambda_n} |\cf_{n+1} f(\tau_\ell \lambda)|^2\\
&= \sum_\ell \sum_\lambda |(\cf_n T_\ell f)(\lambda)|^2
\endalign$$
by Lemma 3.1.4, but then the inductive hypothesis yields
$$
k\left( \frac{d}{|B|}\right)^n \sum_\ell \| T_\ell f\|_{\mu_n}^2
\leq \sum_{\xi\in \Lambda_{n+1}} |\cf_{n+1} f(\xi)|^2
\leq K\left(\frac{D}{|B|}\right)^n \sum_\ell \| T_\ell f\| _{\mu_n}^2.
$$
Combining this with \rom{(*)} gives the desired frame estimate.\qed\enddemo

\subhead 3.2 Limiting estimates \endsubhead
Assume now expansitivity, i.e., that the eigenvalues of $R$ all have modulus $>1$. Then, as noted earlier, the sequence $(\mu_n)$ converges to a certain limit measure $\mu$, so that $\| f\|_{\mu_n}^2$ converges to $\| f\|_\mu^2$ for all continuous $f$. We could then ask for the existence of a {\it renormalization constant\/} $T$ such that the terms
$$
T^n \sum_{\lambda\in \Lambda_n} |\cf_n f(\lambda)|^2\tag 3.13
$$
are comparable to $\| f\|_{\mu_n}^2$ for all $n$.

The next three results present a solution to this circle of asymptotic problems.

Theorem 3.1.5 indicateds that this will be the case for $T=\frac{D}{|B|}$ provided $d=D$; i.e., provided
$$
\langle v, \cm^*\cm v \rangle_{\ell_2^{|B|}} = D\| v\|_{\ell_2^{|B|}}^2 \tag 3.14
$$
for all $v\in \ell_2^{|B|}$. In this case the columns of $\cm$ are orthogonal and $D=|L|$, the last fact used that all the entries in $\cm$ have modulus $1$. In particular, $d=D$ if and only if $\sum_{\ell\in L} \langle b-b',\ell\rangle =0$ for all $b$, $b'\in B$ with $b\neq b'$.

\proclaim{Proposition 3.2.1} Let $(\Omega,\Lambda)$ be a generalized spectral pair, with frame constants $a\leq A$. Suppose \rom{(i)--(v)} are satisfied. If $\sum_{\ell\in L} \langle b-b',\ell \rangle =0$ for all $b$, $b'\in B$ with $b\neq b'$, then
$$
\frac{a}{m_d(\Omega)} \| f\|_\mu ^2 
\,\leq\, \liminf_{n\rightarrow \infty}\left(\frac{|B|}{|L|}\right)^n
\sum_{\lambda\in\Lambda_n} |\cf_n f(\lambda)|^2 \tag 3.15
$$
and
$$
\limsup_{n\rightarrow\infty} \left(\frac{|B|}{|L|}\right)^n
\sum_{\lambda \in \Lambda_n} |\cf_n f(\lambda)|^2
\,\leq\, \frac{A}{m_d (\Omega)} \| f\|_\mu ^2 \tag 3.16
$$
for all $f\in C(\br)$.\endproclaim

\proclaim{Corollary 3.2.2} If the original frame $(\Omega,\Lambda)$ is tight (i.e., $a=A$), then
$$
\lim_{n\rightarrow\infty} \left(\frac{|B|}{|L|}\right)^n 
\sum_{\lambda\in\Lambda_n}|\cf_n f(\lambda)|^2 =\frac{A}{m_d(\Omega)} \| f\|_\mu^2. \tag 3.17
$$
In general $\mu_n \rightarrow \mu$ as $n \rightarrow\infty$, hence $\|f\|_{\mu_n} \rightarrow
\|f\|_\mu$, and $\cf_n f \rightarrow \cf_\mu f$, as $n\rightarrow\infty$ for any continuous function $f$.
\endproclaim

\example{Examples 3.2.3} We will use $(\Omega,\Lambda) = ((0,1),\bz)$ in the examples. Then
$$
\sum_{\lambda\in\bz} |\hat f(\lambda)|^2 =\|f \|_{\cl^2(0,1)}^2
$$
for all $f\in \clo$; hence $a=A=m(\Omega)=1$. Also $K=K^\circ =\bz$; hence the conditions \rom{(i)--(v)} are: $B\subset \bz$, $R\in \bz\backslash(0)$, $L\subset \frac{1}{R}\bz$, if $\ell$, $\ell' \in L$ then $\ell - \ell' \notin \bz$, and the matrix $(\langle b,\ell \rangle)_{b\in B,\ell\in L}$ has rank $|B|$. We will always assume $0\in B$, $0\in L$.

\roster
\item"{(i)}" Let $B=\{0,2\}$, $L=\{0,\frac{1}{4}\}$, $R=4$. Then $\cm = (\langle b,\ell \rangle ) = \pmatrix 1&1\\1&-1\endpmatrix$, hence $d=D=|L|=|B|=2$, so that
$$
\sum_{\lambda\in\Lambda_n} |\cf_n f(\lambda)|^2 \rightarrow \| f\|_\mu^2 \tag 3.18
$$
as $n \rightarrow \infty$, for any continuous $f$. $\Lambda_0\supset \Lambda_1 \supset \Lambda_2 \supset \cdots$, and $\bigcap_{n=0}^\infty \Lambda_n =\cl(L)$, where
$$
\cl(L) =\left\{\sum_{j=1}^p R^{*j} \ell_j : \ell_j \in L,\ p\in\bn \right\}. \tag 3.19
$$
It is known (see \cite{JoPe6}) that
$$
\sum_{\lambda\in \cl(L)} |\cf_\mu f(\lambda)|^2 =0
$$
{\it is\/} possible for $f\neq 0$.

\item"{(ii)}" Same as \rom{(i)}, except we replace $L=\{0,\frac{1}{4}\}$ by $L=\{0,\frac{1}{4}, \frac{2}{4}, \frac{3}{4}\}$. Then $\Lambda_n =\Lambda$ for all $n$. (The sequence $\mu_n$ is the same as the one in \rom{(i)}.) In this example $\sum_{\lambda\in \bz} |\cf_\mu f(\lambda)|^2 = +\infty$ is possible, e.g., if $f(x)=1$ all $x$, then $\cf_\mu f= \hat\mu$ and it is not hard to show that $0\neq \hat\mu(2) = \hat\mu(2\cdot 4^j)$ for any $j\in \bn$.

Now $\cm=\pmatrix 1&1\\1&-1\\1&1\\1&-1 \endpmatrix$, so $d=D=|L|=4$, $|B|=2$, and therefore
$$
\left(\frac{1}{2} \right)^n \sum_{\lambda\in \bz} |\cf_n f(\lambda)|^2 \rightarrow \| f\|_\mu^2 \tag 3.20
$$
as $n\rightarrow \infty$ (for any continuous $f$).

\item"{(iii)}" Let $B=\{0,2\}$, $L=\{0,1/3\}$, $R=3$. The $\cm = \pmatrix 1&1\\ 1&e^{i4\pi/3}\endpmatrix$; hence $d=1$ and $D=3$, and therefore
$$
\left( \frac{1}{2}\right)^n \| f\|_{\mu_n} 
\leq \sum_{\lambda\in\Lambda_n} |(\cf_n f)(\lambda)|^2
\leq \left(\frac{3}{2}\right)^n \| f\|_{\mu_n}.
$$
So in this case, renormalization appears to be impossible. In this example the limit measure $\mu$ is the ordinary Cantor measure on the middle-thirds Cantor-set.

\item"{(iv)}" Same as \rom{(iii)} except $L=\{0,1/3\}$ is replaced by $L=\{0,\frac{1}{3},\frac{2}{3}\}$. (So $\mu_n$ and $\mu$ are the ones from \rom{(iii)}, but $\Lambda_n$ is changed.) Then $\Lambda_n=\Lambda=\bz$ for all $n$. Also $\cm=\pmatrix 1&1\\ 1&e^{i4\pi/3}\\ 1&e^{i8\pi/3} \endpmatrix$; hence $\cm^*\cm =\pmatrix 3&0\\0&3 \endpmatrix$. It follows that $d=D=|L|=3$; hence
$$
\left( \frac{2}{3}\right)^n \sum_{\lambda\in\bz} |\cf_n f(\lambda)|^2
\rightarrow \| f\|_\mu^2 \tag 3.21
$$
as $n\rightarrow\infty$, for all continuous $f$.

\item"{(v)}" Suppose we pick $L$ maximal, i.e., as a complete set of representatives for the cosets $\frac{1}{R}\bz /\bz$. If $b\in B$, $b\neq 0$, then the rank condition prevents $R$ from dividing $b$; that is, it means that $a=e^{i2\pi b/R}\neq 1$. Hence $\sum_{j=0}^{R-1}a^j =\frac{a^R-1}{a-1}=0$, thus $\cm^*\cm =|L|I$. We then get
$$
\lim_{n\rightarrow\infty} \left(\frac{|B|}{|L|}\right)^n \sum_{\lambda\in\bz}
|\cf_n f(\lambda)|^2 = \|f\|_\mu^2 \tag 3.22
$$
for all continuous $f$. (See \rom{(ii)} and \rom{(iv)} above.)\endroster
\endexample

\head Acknowledgments \endhead
The research going into the present paper was done while the coauthors participated in the two work conferences on wavelets and fractals in the Spring of 1994, one at the University of Pittsburgh (USA) in May, and the other at Finsterbergen (near Jena Universit\"at in Germany) in June, and we gratefully acknowledge generous financial support from the respective sponsors of the two conferences. We also thank the individual organizers for stimulating meetings, especially Professors K.-S. Lau (University of Pittsburgh), and C. Bandt (Greifswald Universit\"at). We had helpful discussions with them, with the other participants, and especially with M.~Lapidus, R. Strichartz, D. Mauldin, J. Lagarias, and Y. Wang. We are also grateful to the last mentioned two (Lagarias and Wang) for sending us in July a preprint version of their wonderful paper \cite{LaWa} which inspired and motivated some of the present results.

The work was supported in part by grants from the U.S. National Science Foundation, and from NATO. The first named author was also supported by a University of Iowa Faculty Scholarship.

\enddocument